\documentclass[leqno,11pt]{amsart}
\usepackage{amsmath}
\usepackage{amssymb,amscd}
\usepackage{amsthm}
\usepackage{latexsym}
\usepackage[myheadings]{fullpage}
\usepackage{color}  

\newtheorem{theorem}{Theorem}[section]
\newtheorem{lemma}[theorem]{Lemma}

\newtheorem{proposition}[theorem]{Proposition}
\newtheorem{corollary}[theorem]{Corollary}
\newtheorem{remark}[theorem]{Remark}
\newtheorem{example}[theorem]{Example}

\usepackage{amssymb,graphicx,
	amscd,latexsym,mathrsfs,enumerate,wasysym}

\newcommand{\Ap}{\mathrm{Ap}}
\newcommand{\gr}{\mathrm{gr}}

\title{On the Ratliff-Rush closure of an ideal of a one-dimensional ring}
\author{Veronica Crispin Quinonez}  
\address{Department of Mathematics, Uppsala University, S-751 06, Uppsala, Sweden}
	\email{veronica.crispin@math.uu.se}

\author{Marco D'Anna}
\address{Dipartimento di Matematica e Informatica, Universit\`{a} di Catania, V.le A. Doria, 6, I-95125 Catania, Italy}
\email{marco.danna@unict.it}

\author{Vincenzo Micale}
\address{Dipartimento di Matematica e Informatica, Universit\`{a} di Catania, V.le A. Doria, 6, I-95125 Catania, Italy}
\email{vincenzo.micale@unict.it}

\date{}

\begin{document}


\begin{abstract}
Let $I$ be an ideal in a Noetherian ring $R$ and let $\widetilde{I}$ be its Ratliff-Rush closure. In this paper we study the asymptotic Ratliff-Rush number, i.e. $h(I)=\min\{n\in\mathbb N_+ \mid I^m=\widetilde{I^m}, \ \forall \ m\ge n\}$, in the one-dimensional case. Since $1\le h(I)\le r(I)$, where $r(I)$ is the reduction number of $I$, we look for conditions that determine the extremal values of $h(I)$.
\\
\end{abstract}

\maketitle

\section{Introduction}
Given a Noetherian ring $R$ and a regular ideal $I$ (i.e. $I$ contains a non-zero divisor), the \emph{Ratliff-Rush closure} of $I$ is
defined to be
$$
\widetilde{I}:=\bigcup_{n\ge 1}(I^{n+1}:_RI^n).
$$
In \cite[Theorem 2.1]{RR} it is proved that 
$\widetilde I$ is the largest ideal $J$ of $R$ with the property that $I^n=J^n$, for all $n>>0$; moreover, by \cite[Remark 2.3.2]{RR},
the Ratliff-Rush filtration associated to $I$, i.e. $\{\widetilde {I^n}\}$, asymptotically coincides with the $I$-adic filtration:
$I^n=\widetilde{I^n}$, for every $n>>0$. On the other hand these two filtrations can differ 
greatly in the first steps;
hence controlling these differences can give information on the associated graded ring of $I$. For example, $I^n=\widetilde{I^n}$, for every $n \geq 1$, if and only if $gr_I(R)$ contains a regular element \cite[(1.2)]{HLS}. Ratliff-Rush ideals (i.e. ideals $I$ such that $I=\widetilde I$) and Ratliff-Rush filtrations have been widely studied in the last thirty years in different contexts and from different points of view (see e.g 
\cite{HLS}, \cite{RV}, \cite{RS}, \cite{C} and \cite{AA2}).

In view of the fact that, for any $n$ large enough, $I^n =\widetilde{I^n}$, it is natural to define the \emph{asymptotic Ratliff-Rush number} of $I$, as 
$$
h=h(I)=\min\{n\in\mathbb N_+ \mid I^m=\widetilde{I^m}, \ \forall \ m\ge n\}
$$
($h(I)$ was introduced in \cite{D'A-G-H} as the minimum non
negative integer with the prescribed property,
but, since the equality $I^0=R=\widetilde{I^0}$ always trivially holds, we prefer to restrict to $\mathbb N_+$).

The starting point of our investigation 
is a result in \cite{D'A-G-H}, which states that, if $I$ has a principal reduction $x$, then 
$h(I)\leq r(I)$, where $r(I)$ is the \emph{reduction number} of $I$, i.e. 
$r(I)=\min\left\{n\in\mathbb N\ |\ I^{n+1}=xI^n \right\}$ (it is well known that this number is independent of the principal reduction $x$). It is worth noticing that, for $1\leq m < n < r(I)$, it can happen that $I^m=\widetilde{I^m}$ and $I^n \neq \widetilde{I^n}$ (see, e.g. Example \ref{h=r}). 
The natural context for which every regular ideal has a principal reduction is the one-dimensional case. Assuming also that $R$ is local, with infinite residue field, we get that the minimal reductions of $I$ coincide with the principal ones. We also notice that it would be interesting to understand if, in the general (non-local) one-dimensional Noetherian case, there is a relation between $h(I)$ and $r(I)$; this is not the case in higher dimension as shown in \cite[Remark 2.7]{AA1}.

So, from now on, we will assume that $R$
is a one-dimensional, local, Noetherian ring 
with regular maximal ideal $\mathfrak m$
(that, in this case, means $R$ Cohen-Macaulay)
and infinite residue field,
and we will look for conditions on $I$ that
imply or characterize the extremal values of $h(I)$, i.e. either $h(I)=1$ or $h(I)=r(I)$.

In order to study the case $h(I)=1$, that in our setting, by \cite[(1.2)]{HLS}, is equivalent to say that  $gr_I(R)$  is Cohen-Macaulay, we make use of a pullback construction;
more precisely, we assume that the residue field $k=R/\mathfrak m \subseteq R$ and we consider the new ring $U=k+I$, which is again 
a one-dimensional, local, Noetherian ring, with maximal ideal $\mathfrak n=I$, that is also a regular ideal.
Multypling $R$ by the principal reduction $x$
of $I$, we obtain the ideal $J=xR$ and it turns out that $\gr_I(R) \cong \gr_{\mathfrak n}(J)$, both as $\gr_{\mathfrak n}(U)$-module and as a $\gr_I(R)$-module. Hence we can read the condition $h(I)=1$, checking if 
$\gr_{\mathfrak n}(J)$ is a Cohen Macaulay $\gr_{\mathfrak n}(U)$-module (see Corollary \ref{h=1}). This checking can be done with an effective computation in the case of monomial ideals of a numerical semigroup ring, using a result in \cite{D'A-J-St} (see Theorem \ref{comp}).
It remains open to understand if this idea can be generalized to a larger class (e.g. to analytically irreducible domains). Finally, using our Theorem \ref{comp} we are able to check that 
we can have $h(I)=1$, with $r(I)$ arbitrarily large (see Example \ref{ex}).

As for the case $h(I)=r(I)$, we turn back to the more general one-dimensional context, assuming that the integral closure, $\overline R$ of $R$,  
in its total ring of fractions $Q(R)$ is
a finite $R$-module. This fact implies that 
the conductor ideal $C:=(R:_{Q(R)}\overline R)$ is a regular 
ideal of $R$. Under these hypotheses we deepen the study initiated in \cite{D'A-G-H}, giving a
new sufficient condition that implies $h(I)=r(I)$ (see Proposition \ref{1}). Finally, we specialize again our results to the case of numerical semigroup rings, obtaining a new numerical sufficient condition (see Proposition \ref{suff}). 


\section{Preliminaries on numerical semigroup rings}

In this section we collect some basic notions and results on numerical semigroups and numerical semigroup rings that we will use in the sequel. For the proofs of the stated results we refer to \cite{A-D-GS}.

A \emph{numerical semigroup} $S$ is a submonoid of $(\mathbb N,+)$ such that $|\mathbb N\setminus S|$ is finite; the smallest integer $c$ such that $c+\mathbb N \subseteq S$ is called  the \emph{conductor} of $S$. It is well known that $S$ is finitely generated and has a unique minimal system of generators. Throughout the whole paper, $S=\left\langle n_1,\dots,n_\nu\right\rangle$ is a numerical semigroup minimally generated by $n_1<\cdots<n_\nu$; the ring $k[[S]]=k[[t^{n_1}, \dots, t^{n_\nu}]]$
is the corresponding numerical semigroup ring, with maximal ideal $\mathfrak{m}=(t^{n_1}, \dots, t^{n_\nu})$. The smallest nonzero element of
$S$, $n_1$, is called the \emph{multiplicity of} $S$ and is denoted by $m$; it is well known that $m=e(k[[S]])$, the multiplicity of $k[[S]]$. 

A 
\emph{relative ideal} of $S$ is a non-empty set $E$ of integers such that $E+S\subseteq E$ and $s+E\subseteq S$ for some $s\in S$; when it is contained in $S$, $E$ is simply called an \emph{ideal} of $S$. 
As for semigroups, we define the multiplicity of a relative ideal $E$ as the smallest element of $E$ and we denote it with $e(E)$. Note that for relative ideals $E_1$ and $E_2$ of $S$, the set $E_1+E_2=\{e_1+e_2\ |\ e_1\in E_1, e_2\in E_2\}$ is also a relative ideal. In particular, for $z\in\mathbb Z$, $z+S=\{z+s\ |\ s\in S\}$ is the principal relative ideal of $S$ generated by $z$. For any ideal $E$ of $S$, we can always express it as $E=(e_1+S)\cup\cdots\cup(e_h+S)$, for some $e_i\in E$; then, we write $E=\{e_1,\dots,e_h\}+S$ and we can always assume that the set $\{e_1,\dots,e_h\}$ is minimal, i.e., for all $i=1,\dots,h$, $e_i\notin \bigcup_{j\ne i}(e_j+S)$; it is straightforward to see that $E$ has a unique minimal set of generators. Moreover we always enumerate the generators in increasing order;
so, in particular, $e(E)=e_1$. By difference of two ideals $E_1$ and $E_2$, we mean the ideal $E_1-E_2=\{z\in\mathbb Z\ |\ z+E_2\subseteq E_1\}$. 

We denote by $M=S\setminus\{0\}$ the maximal ideal of $S$ and we set
$lM=M+\cdots+M$. The \emph{blow-up} of $S$ is defined as the numerical semigroup $B(S)=\bigcup_l(lM-lM)=\left\langle m, n_2-m,\dots, n_\nu-m\right\rangle$. It is well known that $B(S)=lM-lM=lM-lm$ for $l$ large enough.

Let $\omega_i=\min\{s\in S\ |\ s\equiv i\ (mod\ m)\}$. The \emph{Ap\'{e}ry set} of $S$ with respect to $m$ is the set $\Ap_m(S)=\{\omega_0=0, \omega_1,\dots, \omega_{m-1}\}$. In the same way we denote $\Ap_m(B(S))=\{\omega'_0=0, \omega'_1,\dots, \omega'_{m-1}\}$. It follows from the definition that $\omega_i\ge\omega'_i$ for all $i=0,\dots, m-1$ and we define the microinvariants of $S$ as the integers $a_i(S)$ such that $\omega'_i+ma_i(S)=\omega_i$. Moreover, we set $b_i(S)=\max\{l\ |\ \omega_i\in lM\}$. A criterion for the Cohen-Macaulayness of the associated graded ring, proved by Barucci and Fr\"oberg for analytically irreducible domains, implies the following result.

\begin{theorem} \cite[Theorem 2.6]{BF}
	The associated graded ring $\gr_{\mathfrak m}(k[[S]])$ is Cohen-Macaulay if and only if $a_i(S)=b_i(S)$ for each $i=0,\dots, m-1$.
\end{theorem}

Let $E=\left\{e_1,\dots,e_n\right\}+S$ be an ideal of a semigroup $S$. 
The \emph{Ap\'{e}ry set} of $E$, with respect to the multiplicity  $m$ of $S$, is $\Ap_m(E)=\{\alpha_0,\alpha_1,\dots, \alpha_{m-1}\}$, where $\alpha_i$ is the smallest element in $E$ congruent to $i$ modulo $m$; we notice that $m$ may not be in $E$. In \cite{D'A-J-St}, the authors define the \emph{blow-up} of $E$ as 
$$
B(E)=\bigcup_{i\ge 1}\left(E+\left(i-1\right)M\right)-iM.
$$
and prove that $B(E)=\left(E+\left(i-1\right)M\right) -iM$ for $i$ large enough.

\begin{lemma}\label{B} \cite[Lemma 3.1]{D'A-J-St}
	Let $S$ be a numerical semigroup with maximal ideal $M$ and multiplicity $m$. Then $B(E)=\left\{e_1-m,\dots,e_n-m\right\}+B(S)$.
\end{lemma}

\begin{remark}
	If $m\in E$, then $0\in B(E)$ and so $B(E)=B(S)$.
\end{remark}

Let $\Ap_m(B(E))=\{\alpha'_0,\alpha'_1,\dots, \alpha'_{m-1}\}$. As for the semigroup case, we define the microinvariants of $E$ as the integers $a_i(E)$ such that $\alpha'_i+ma_i(E)=\alpha_i$. Moreover, we set $b_i(E)=\max\{l+1\ |\ \alpha_i\in lM+E\}$. It easy to see that $a_i(E)\ge b_i(E)$ \cite[Remark 3.2]{D'A-J-St}.

With these notations it is possible to generalize Theorem 1 for ideals.
\begin{proposition} \label{D-J-S}\cite[Proposition 3.6]{D'A-J-St}\label{ECM} Let $E=\{e_1, \dots, e_n\}+S$ and let $I=(t^{e_1},	\dots , t^{e_n})$.
	The following statements are equivalent:
	\begin{enumerate}
		\item  $\gr_{\mathfrak m}(I)$ is a one-dimensional  Cohen-Macaulay $\gr_{\mathfrak m}(k[[S]])$-module;
		\item $a_i(E)=b_i(E)$ for all $i=0, \ldots, m-1$.
	\end{enumerate}
\end{proposition}

\section{The case $h(I)=1$}

Let $(R,\mathfrak m)$ be a one-dimensional, Noetherian, Cohen-Macaulay local ring, with residue field $k=R/\mathfrak m\subseteq R$. Let $I$ be a regular ideal of $R$ (i.e. $I$ contains a nonzero divisor) and let $x$ be a minimal reduction of $I$. Let us consider the 
subring $k+I$ of $R$,
that can be viewed as a pullback as shown 
by the following commutative square
(where $k+I=\pi^{-1}(k)$):

$$  \hskip1cm
\CD k+I  @>>> k \\
@VVV   @VVV \\
R  @>\pi>> R/I 
\endCD 
$$
Since $R/I$ is a finite $k$-vector space,
being an artinian ring, also the inclusion 
$k+I \subseteq R$ is finite and $k+I$ is Noetherian (see \cite[Proposition 1.8]{fo}); from this fact, it also follows immediately that $U=k+I$ is one-dimensional and local (with maximal ideal $\mathfrak n=I$). Moreover, being $I$ a regular ideal, $U$ is Cohen-Macaulay.

Multiplying $R$ by $x$ we obtain an ideal of $U$ that, from now on, we denote by $J:=xR$.

\begin{proposition} Under the hypotheses and notations introduced above,
	$\gr_I(R)$ is isomorphic to $\gr_{\mathfrak n}(J)$, both as a $\gr_{\mathfrak n}(U)$-module and as a $\gr_I(R)$-module.
\end{proposition}
\begin{proof}
	We note that $\gr_I(R)=\frac{R}{I}\oplus	\frac{I}{I^2}\oplus\cdots$ and, by $\mathfrak n=I$, we have that 
	$$
	\gr_{\mathfrak n}(J) \ = \ \frac{J}{\mathfrak nJ}\oplus\frac{\mathfrak nJ}{\mathfrak n^2J}\oplus\cdots \ = \ \frac{xR}{xI}\oplus\frac{xI}{xI^2}\oplus\cdots.
	$$ 
	%
	As $x$ is a regular element in both $R$ and $U$, we get that $I^{h}\cong xI^{h}$, for any $h \in \mathbb N$, as $U$- and $R$-module. Hence we get $\frac{I^j}{I^{j+1}}\cong\frac{xI^j}{xI^{j+1}}$ as $U$- and $R$-module. Thus $\gr_{\mathfrak n}(J)$ is isomorphic to $\gr_I(R)$, both as a $\gr_{\mathfrak n}(U)$-module and as a 
	$\gr_I(R)$-module.
\end{proof}

In light of the above isomorphism we immediately obtain
the following result.

\begin{corollary} Preserving the hypotheses and notations of the previous proposition,
	$\gr_I(R)$ is a Cohen-Macaulay ring if and only if $\gr_{\mathfrak n}(J)$ is a Cohen-Macaulay $\gr_{\mathfrak n}(U)$-module. 
\end{corollary}


Remembering that in our setting the condition $\gr_I(R)$ Cohen-Macaulay is equivalent to say $h(I)=1$, we can rephrase the previous corollary.

\begin{corollary}\label{h=1} Under the standing hypotheses and notations, $h(I)=1$ if and only if $\gr_{\mathfrak n}(J)$ is a Cohen-Macaulay $\gr_{\mathfrak n}(U)$-module.
\end{corollary} 

If we restrict to the case of numerical semigroup rings, the previous corollary, together with Proposition \ref{D-J-S}, 
produces an effective computational method to
check when $h(I)=1$.

More precisely, let us fix the following notations through the rest of this section:
let $S$ be a numerical semigroup and let $E$ be an ideal of $S$ with multiplicity $e=e(E)$. Let us consider the numerical semigroup $T=\{0\}\cup E$ (hence the maximal ideal of $T$ coincides with $E$) and let $F=e+S\subseteq T$. Is easy to see that $F$ is an ideal of $T$ with $e(F)=e$.

Let $R=k[[S]]$ and $I=(t^a\ |\ a\in E)$. Set $x=t^{e}$, $U=k[[T]]$; clearly $U=k+I$ so we are in the setting of the beginning of this section. Finally, set $J=xR$, which is an ideal of $U$. Under these assumptions,
if $v$ is the usual discrete valuation on $k[[t]]$, we get: $v(R)=S$, $v(I)=E$, $v(U)=T$ and $v(J)=F$. 

Let $\left\{f_1,f_2,\dots,f_{n}\right\}$ be the generators of $F$ as ideal in $T$ (notice that $f_1=e$) and let $B(F)$ and $B(T)$ be the blow-up of $F$ and $T$, respectively.

\begin{remark} 
	Since $f_1=e$, Lemma \ref{B} implies that the blow-up of $F$ as ideal of the semigroup $T$ is 
	$$
	B(F)=\left\{ f_1-e,f_2-e,\dots, f_{n}-e\right\}+B(T)=B(T)
	$$
\end{remark}

The next result gives the promised computational method that allows to check when $h(I)=1$.

\begin{theorem}\label{comp} 
	Under the standing hypotheses and notations, $h(I)=1$ if and only if $a_i(F)=b_i(F)$, for every $i=0,\dots,e-1$.
\end{theorem}
\begin{proof} 
	By Proposition \ref{D-J-S}, we have that  $a_i(F)=b_i(F)$, for every $i=0,\dots,e-1$, if and only if $\gr_{\mathfrak n}(J)$ is a Cohen-Macaulay 
	$\gr_{\mathfrak n}(U)$-module. By Corollary
	\ref{h=1} we immediately obtain the thesis.
\end{proof}

\begin{example}\label{ex}
	Set $I=(t^9,t^{11}) \subset k[[S]]$ with $S=\left\langle 6,9,11\right\rangle$. Let us show that $h(I)=1$ using last theorem. 
	
	We have $v(I)=E=\{9,11\}+S=\{9,11,15,17,18,20, 21,22, 23,24,26,\rightarrow\}$ (where the arrow means that all the integer greater than $26$ belongs to $E$), and $3E=9+2E$; from this fact it follows immediately that $r(I)=2$. Moreover $T=\{0\}\cup E$
	and $$F=9+S=\{9,15,18,20,21,24,26,27,29,30,31,32,33,35,\rightarrow\};$$ thus $\Ap_9(F)=\{9,37,20,21,31,32,15,43,26\}$. 
	
	Furthermore, by $T=\left\langle 9,11,15,17,21,23\right\rangle$, we easily get  $B(F)=B(T)=\{0,2,4,6,8,\rightarrow\}$ and $\Ap_9(B(F))=\{0,10,2,12,4,14,6,16,8\}$; hence $a_0(F)=a_3(F)=a_6(F)=1$, $a_2(F)=a_5(F)=a_8(F)=2$ and $a_1(F)=a_4(F)=a_7(F)=3$. 
	
	The maximal ideal of $T$ is $E$ and so $$E+F=\{18,20,24,26,27,29,30,31,32,33,\rightarrow\}$$ and $lE+F=\{(l+1)\cdot 9, (l+1)\cdot 9+2, (l+1)\cdot 9+4,(l+1)\cdot 9+6,(l+1)\cdot 9+8,\rightarrow\}$ for every $l\ge 2$. This implies that $\alpha_i\notin E+F$ for $i\in\{0,3,6\}$, $\alpha_i\in (E+F)\setminus (2E+F)$ for $i=2,5,8$ and $\alpha_i\in (2E+F)\setminus (3E+F)$ for $i=1,4,7$. Therefore, we get $a_i(F)=b_i(F)$ for $i=0,\dots,8$.
\end{example}

In the previous example one could compute 
$h(I)$ directly checking that $I=\widetilde I$,
since $r(I)=2$ and so all the subsequent powers of $I$ have to be Ratliff-Rush closed.
Hence it is clear that the computational method given by Theorem \ref{comp} becomes convenient when the reduction number of $I$ increases. 

In the following example we consider a family of semigroups $S_n$ (with $n \geq 3$) and ideals $E_n \subset S_n$, introduced in
\cite[Example 2.3]{D'A-G-H2}, that have the following properties: $S_n$ is minimally generated by $n+1$ elements; $hE_n$ is minimally $n$-generated for every $h \geq 1$; the reduction number of $E_n$ is $n-1$. 
We will show that,  setting 
$I_n=(t^x: x \in E_n) \subset k[[S_n]]$ and using Theorem \ref{comp}, it is possible to prove that $h(I_n)=1$, for any $n \geq 3$.
Hence, as byproduct, we obtain that the asymptotic Ratliff-Rush number of an ideal,  with reduction number arbitrarily large, can be equal to $1$. Notice that in higher dimension such a situation has been proved in \cite[Example 3.2]{AA2}.

\begin{example}\label{ex}
	Fix $n\geq 3$ and set $S_n = \langle a,b,d,c_3,\dots , c_n\rangle$, where $a=2n$,
	$b=4n-1$, $d=n(2n-1)$ and $c_h=(n+h)(2n-1)+1$, for any $h=3,\dots n$. 
	
	Let $E_n=\{a,b,c_3,\dots, c_n\}+S_n$ and set
	$I_n=(t^x: x \in E_n) \subset k[[S_n]]$.
	We will show that 
	$h(I_n)=1$ using Theorem \ref{comp}. 
	
	We have: 
	$$
	T_n=\{0\}\cup E_n =\langle a,b,c_3,\dots, c_n,a+d,b+d, c_3+d, \dots, c_n+d\rangle
	$$
	but, since $c_h+d=(h-1)b+(2n+1-h)a$ for any $h=3,\dots n$, the last generators are superfluous. Hence 
	$T_n=\langle a,b,c_3,\dots, c_n,a+d,b+d\rangle$.
	Since $d$ is the only generator of $S_n$ not belonging to $E_n$, we easily obtain that  $F_n=S_n+a$ is generated, as ideal of $T_n$, by $\{a, a+d\}$. 
	
	As proved in \cite[Example 2.3]{D'A-G-H2}, $hE_n$ is minimally generated, as ideal of $S_n$, by $\{ha, (h-1)a+b, \dots, hb, (h-1)a+c_{h+2}, \dots, (h-1)a+c_n\}$, for any $h < n-1$, while $(n-1)E_n=\{(n-1)a, (n-2)a+b, \dots, (n-1)b\}+S_n$ and $nE_n=(n-1)E_n+a$.
	
	Using induction and the relations 
	between the generators proved in \cite[Example 2.3]{D'A-G-H2}, it is possible to show that $hb-a \notin S_n$ (and so $hb-a \notin T_n$) for every $h=1,\dots n-1$; therefore the Ap\'ery set of $T_n$ (ordering the elements by their residue modulo $a$) is the following:
	$$
	\Ap_a(T_n)=\{0,c_n, c_{n-1},\dots, c_3, d+b,d+a,(n-1)b, \dots, 2b,b\}.
	$$
	From this fact one can show that the Ap\'ery set of $F_n=a+S_n$ is 
	$$
	\Ap_a(F_n)=\{a,a+c_n, a+c_{n-1},\dots, a+c_3, a+d+b,d+a,a+(n-1)b, \dots, a+2b,a+b\}.
	$$
	Finally, since $d=n(b-a)$, we get $B(F_n)=B(T_n)=\langle a,b-a, c_3-a, \dots, c_n-a\rangle$, whose Ap\'ery set with respect to $a$ is 
	$$
	\Ap_a(B(F_n))=\{0,c_n-a, c_{n-1}-a,\dots, c_3-a,
	 d+b-a, d,$$
	 $$(n-1)(b-a), \dots, 2(b-a),b-a\}
	$$
	and, using the computation of the $hE_n$, that are the multiples of the maximal ideal of $T_n$,
	one can show that the orders $b_i$ of the elements of the Ap\'ery set of $F_n$ are
	$$
	b_0=0, b_1=\dots=b_{n-1}=2, b_n=1, b_{n+1}=n, b_{n+2}=n-1, \dots , b_{2n-1}=2,
	$$
	that, therefore, coincide with the $a_i$.
\end{example}

%
%
%
%
%

If $I$ is integrally closed, we can characterize $h(I)=1$ in terms of $\gr_{\mathfrak n}(U)$ (instead of in terms of
$\gr_{\mathfrak n}(J)$).

\begin{proposition} We preserve the hypotheses and notations introduced before Theorem \ref{comp}.
	Assume that $I$ is an integrally closed
	ideal of $k[[S]]$. Then $h(I)=1$ if and only if $\gr_{\mathfrak n}(U)$ is Cohen-Macaulay.
\end{proposition}
\begin{proof}
	Since $I$ is integrally closed, we have that $E=\overline{E}$ where $\overline{E}$ is the integral closure of $E$, that is $\overline{E}=\left\{s\in S\ |\ s\ge e\right\}$. Hence, $e+s\in E\setminus 2E$, for every $s\in S\setminus E$, as every such $s$ is smaller than or equal to $e$. This implies that $t^e$ is a non zero divisor for $\gr_I(R)$ if and only if $t^e$ is a non zero divisor for $\gr_{\mathfrak n}(U)$.   
\end{proof}

It is well known that $I=\overline{I}$ implies $I=\widetilde{I}$. There are examples, in a more general context, that show that the inverse does not hold. The inverse does not hold also in the numerical semigroup context even if $I$ is contained in the conductor of $k[[S]]$.

\begin{example}
	Set $I=(t^9,t^{11}) \subset k[[S]]$ with $S=\left\langle 4,5,6\right\rangle$. We have  
	$v(I)=E=\{9,11\}+S=\left\{9,11,13,\stackrel{}{\rightarrow}\right\}$; therefore $2E=9+E$, so $r(I)=1$ and, thus, also $h(I)=1$. However, $E$ is included in the conductor of $S$ and $\overline{E}\setminus E=\left\{10,12\right\}$. 
\end{example} 

We conclude this section with a consequence of our construction that holds in the general case.

\begin{remark} Under the hypotheses and notations of the beginning of this section,
	since $U=k+I$ is a local Noetherian one-dimensional ring, its embedding dimension
	$\nu(U)$ is bounded above by its multiplicity
	$e(U)$.
	It is well known that being of maximal embedding dimension (i.e. $\nu(U)=e(U)$)
	is equivalent to the stability of the maximal ideal, that, in our case, is $I$.
	Hence, if $U$ is a ring of maximal embedding dimension, we have that $I^2=xI$, and therefore 
	it is straightforward to check that the image of $x$ in $\gr_I(R)$ is a non-zerodivisor, i.e. $\gr_I(R)$ is a Cohen Macaulay ring, that is $h(I)=1$.
\end{remark}

\section{The case $h(I)=r(I)$}

Let $R$ be a one-dimensional, local, reduced Noetherian ring having total ring of fractions $Q(R)$ and let assume that the integral closure $\overline {R}$ of $R$ in $Q(R)$ is a finitely generated $R$-module. Let $C=(R:_{Q(R)}\overline{R})$ denote the conductor of $\overline{R}$ into $R$.

Let $x$ be a principal reduction of $I$ and let $r=r(I)=\min\{n\in\mathbb N\ |\ I^{n+1}=xI^n\}$ be the reduction number of $I$. As noticed in \cite[Discussion 2.5]{D'A-G-H}, if $I$ is a regular ideal with a principal reduction $xR$, then, for every $s\ge r$,
$$
\widetilde{I}=(I^{s+1}:I^s)=(xI^s:_{Q(R)}I^s)\cap R=xR^I\cap R=I^sx^{-s+1}\cap R,
$$
where $R^I$ denotes the blow-up of $I$, i.e.,
in our setting, $R^I=\bigcup_{n\in \mathbb N}(I^n:I^n)$; it is well known that $R^I\subseteq \overline R$ and that $R^I=(I^s:I^s)$, for any $s \geq r(I)$ (see again \cite[Discussion 2.5]{D'A-G-H}).

From the equality $I^{r+1}=xI^r$ it follows immediately that
$I^{r+n}=I^rx^n$ for every $n\in\mathbb N$.
Moreover, it is straightforward to check that $x^m$ is a principal reduction of $I^m$ for every $m\in\mathbb N$.
Using this facts, we can 
obtain a bound on the reduction number of $I^m$.


\begin{lemma} 
	Fix and integer $m\in\mathbb N$, $m\geq 2$ and let $l\in\left\{0,\dots,m-1\right\}$ be such that $m$ divide $r+l$. Then $r(I^m)\le\frac{r+l}{m}$.
\end{lemma}
\begin{proof}
	The thesis follows by the following chain of equalities: 	
	$$ \left(I^{m}\right)^{\frac{r+l}{m}+1}=I^{r+l}I^{m}=I^mI^rx^{l}=I^{r+m}x^{l}=
	I^lI^{r+(m-l)}x^{l}= 
	$$
	$$
	\ \ \ =I^lI^rx^{m-l}x^l=I^{r+l}x^m=\left(I^{m}\right)^{\frac{r+l}{m}}x^m .
$$
\end{proof}

\begin{lemma}
	Let $m\in\mathbb N$, $m \geq 1$; then $I^rx^{-r+m}\cap R=\widetilde{I^m}$.
\end{lemma}
\begin{proof}
	If $m=1$ the thesis is given by the formula
	at the beginning of the section. If $m \geq 2$,
	let $l\in\left\{0,\dots,m-1\right\}$ such that $m$ divide $r+l$. Then, making use of the previous lemma, we get 
	$$ I^rx^{-r+m}\cap R=I^rx^lx^{-r+(m-l)}\cap R=I^{r+l}x^{-r+(m-l)}\cap R= $$
	$$=I^{r+l}x^{-(r+l)+m}\cap R=\left(I^{m}\right)^{\frac{r+l}{m}}\left(x^m\right)^{-\frac{r+l}{m}+1}\cap R=\widetilde{I^m} .$$
\end{proof}

Notice that if an ideal $J$ with principal reduction $yR$ is included in the conductor, by $R^J\subset \overline R$, it follows that $\widetilde J= yR^J\cap R=yR^J$, so the intersection with $R$ is superfluous. 
As observed also in \cite{D'A-G-H}, this means
that an ideal included in the conductor
is Ratliff-Rush closed if and only if it is stable (i.e. it has reduction number $1$).
If we apply this remark to $I^m$ we obtain the following result.

\begin{corollary}\label{l<r} Let $I$ be an ideal with reduction number $r$.
	If there exists $m<r$, such that $I^rx^{-r+m}\subseteq C$, then $I^m\ne \widetilde{I^m}$
\end{corollary}

\begin{proof}
	If $m<r$, then  $I^mx^{r-m}\subsetneq I^r$; therefore $I^m \subsetneq I^r x^{-r+m}$. By the previous lemma and by the hypothesis $I^rx^{-r+m}\subseteq C$, it follows that $I^r x^{-r+m}=\widetilde{I^m}$,
	that, in turn, implies the thesis.
\end{proof}

We can use the previous corollary to improve \cite[Proposition 3.10]{D'A-G-H}.

\begin{proposition}\label{1} Let $I$ be an ideal with reduction number $r$ and 
	set $l=\min\{m\in\mathbb N\ :\ I^rx^{-r+m}\subseteq C\}$. If $l<r$, then $h(I)=r$.
\end{proposition}
\begin{proof}
	By Corollary \ref{l<r}, for any $m \in \{l, \dots , r-1\}$, $I^m\neq \widetilde{I^m}$. The thesis follows immediately.
\end{proof}

\begin{example}\label{h=r}
	Set $I=(t^4,t^5,t^{11})\subset k[[S]]$, with $S=\langle 4,5,11\rangle=\{0,4,5,8,\rightarrow\}$. Clearly $I$ is the maximal ideal of $k[[S]]$, so $I= \widetilde I$.
	On the other hand, $I^2=(t^8,t^9,t^{10})$, $I^3=(t^{12},t^{13},t^{14},t^{15})$ and, for any $k\geq 1$,
	$I^{3+k}=x^{k}I^3$ (where $x=t^4$); thus $r(I)=3$. But 
	$I^2 \subseteq C$, hence $I^2\subsetneq \widetilde{I^2}$; in fact, arguing as in Corollary \ref{l<r}, $\widetilde{I^2}=I^3x^{-1}$, that contains $t^{11}$. 
\end{example}
We conclude the paper applying the above proposition to the case of numerical semigroup rings. So
let $S$ be a numerical semigroup with conductor $c$ and let $R=k[[S]]$;
it is well known that $c=\min\left\{v(x)\ |\ x\in C\right\}$, where $C$ is the conductor of $R$.

\begin{proposition}\label{suff}
	Let $S$ be a numerical semigroup with conductor $c$ and let $E$ be an ideal of $S$ with multiplicity $e=e(E)$. Let $R=k[[S]]$ and let $I=(t^a\ |\ a\in E)$. Assume that $r(I)=r\ge 2$; if $(r-1)e\ge c$, then $h(I)=r$.
\end{proposition}
\begin{proof}
	By assumption, since $x=t^e$ is a principal reduction of $I$, we have $I^rx^{-r+(r-1)}\subseteq C$; therefore the integer $l$ defined in Proposition \ref{1} is such that $l\le r-1<r$. Again by Proposition \ref{1}, we get the thesis.
\end{proof}

\begin{example}
	Set $I=(t^7,t^8) \subset k[[S]]$ where $S=\left\langle 4,5,7\right\rangle=\left\{0,4,5,7,\stackrel{}{\rightarrow}\right\}$; in this case the conductor of $S$ is $c=7$.  Since $E=\{7,8\}+S=\left\{7,8,12,\stackrel{}{\rightarrow}\right\}$, we have $e=7$ and, by simple calculations, $r=4$. Finally, Proposition \ref{suff} implies that $h(I)=4$. 
\end{example}

\end{document}